\documentclass[12pt]{article}
\usepackage{amsmath}
\usepackage{graphicx}
\usepackage{amsthm}
\usepackage{nicefrac}
\usepackage{breakcites}
\usepackage{amssymb}
\usepackage{natbib}
\usepackage{setspace}
\usepackage{color}
\usepackage{tikz}
\usepackage{datatool}
\usepackage{multirow}
\def\normal{\mathcal{N}}
\def\natus{\mathbb{N}}
\def\reals{\mathbb{R}}
\def\E{\text{E}}
\def\Var{\text{Var}}
\def\Cov{\text{Cov}}
\newtheorem{remark}{Remark}
\newtheorem{proposition}{Proposition}

\allowdisplaybreaks

\title{On mean decomposition for summarizing conditional distributions}
 \author{Celia Garc\'{i}a-Pareja$^{1\dagger}$, Matteo Bottai$^{1}$ \\
\\
\normalsize{$^{1}$Unit of Biostatistics, IMM, Karolinska Institutet, Stockholm, Sweden.}\\
\\
\normalsize{$^\dagger$corresponding author: Celia Garc\'{i}a-Pareja, Unit of Biostatistics, IMM, Karolinska Institutet,}\\
\normalsize{Box 210, 171 77 Stockholm, Sweden. E-mail: celia.garcia.pareja@ki.se }
}

\date{}

\begin{document} 

\maketitle 
\begin{abstract}
  
We propose a summary measure defined as the expected value of a random variable 
over disjoint subsets of its support that are specified by a given grid of proportions,  
and consider its use in a regression modeling framework. 
The obtained regression coefficients
provide information about the effect of a set of given covariates on the variable's expectation in each specified 
subset. 
We derive asymptotic properties for a general estimation approach that are based on those of 
the chosen quantile function estimator for the underlying probability distribution. A bound on the variance 
of this general estimator is also provided, which relates its precision to the given grid of proportions and that 
of the quantile function estimator, as shown in a simulation example. 
We illustrate the use of our method and its advantages in two real data applications, where we show its potential 
for solving resource-allocation and intervention-evaluation problems.

\end{abstract}
\small{\textit{keywords:} asymptotics; compound expectation; conditional quantile function; regression model; summary measure}
\section{Introduction}
  Statistical summary measures aim to comprise relevant 
information from the distribution of a variable of interest. Despite their attractive simplicity, 
central tendency measures are not always appropriate when we intend to describe the entire distribution. 
The sample mean, for example, provides practical but often insufficient information.
Conversely, a set of sample quantiles can provide a more detailed picture but lack information on how 
the distribution behaves between elements of the set.

Based on the same mathematical relation between these two extremes exploited in \citet{Wang2010}, we propose to summarize the variable's distribution by 
specifying a grid of proportions that divide its mean into a sum of components. From these components, one can easily derive 
the variable's expected value on different segments of its support. Our approach results in a set--valued summary 
measure, which we refer to as compound expectation, that describes the variable's entire distribution in terms of 
expected values, and whose elements relate to specific fractions of the variable's support.

Following this same formulation, the compound expectation can be easily extended to a regression framework 
in which we summarize the variable's conditional distribution given covariates. Similarly to the univariate case, 
we obtain regression coefficients that measure average differences on the outcome variable over distinct fractions 
along its entire distribution. 
Other regression models that characterize entire conditional distributions 
have been previously proposed, e.g., by  
\citet{Koenker1978}, 
\citet{Newey1987} or \citet{Breckling1988}, and their advantages discussed at length 
(\citet{Koenker2005}, \citet{Chambers2006}, \citet{Ehm2016}). 
However, conditional averages are still preferred in numerous applications, where measuring average behavior 
is essential for drawing meaningful conclusions. In this regard, our proposed regression model allows 
detecting different associations between outcome variable and covariates along the distribution, 
while preserving the appealing advantages of the mean.

All the components, which can take on positive or negative values, sum to the total expectation, and one can identify their contributions to it. Such contributions are 
measured as the proportion each component represents of the total mean. 
This feature is of convenience in numerous settings, 
where interests focus on how specific fractions of a 
variable influence its average. Examples 
include how premature adult-age deaths affect life expectancy \citep{Seaman2016}, 
how extreme warmth raises average temperature \citep{Meehl2000}, 
or how the most skilled increase average intelligence quotient \citep{Wai2011}. 

The compound expectation can also be considered as a measure of concentration, 
where all components being equal indicates no concentration, 
and all but one components being zero indicates maximal concentration \citep{Egghe1990}. 
In such context, one can analyze how the variable accumulates around its different fractions, which is of use for 
solving resource-allocation problems.
Approaches based on concentration curves, which originally appeared as a metric to study the distribution of 
wealth \citep{Lorenz1905}, have been recently used to address these problems,  
as exemplified by their application to transport logistics or epidemiology (\citet{Delbosc2011}, \citet{Mauguen2016}, 
\citet{Christopoulos2017}). Graphical tools, however, narrow modeling possibilities and 
require restrictions or ad hoc ordering criteria to be compared (see \citet{Shorrocks1983}, \citet{Davies1995} and 
\citet{Aaberge2009}, among others). Both these limitations are effectively overcome by the regression 
approach we propose in this paper. 

The rest of this paper is structured as follows. In Section \ref{comp_exp}, we define the compound expectation and related quantities. 
In Section \ref{reg_framework}, 
we extend its use to a regression framework, which serves as a modeling tool in the presence of covariates. 
Our proposed estimation approach and its asymptotic properties are described in 
Section \ref{asymp_properties}. In Section \ref{grid} we present a simulated data example that illustrates how different grids 
of proportions might reflect in final 
results, and in Section \ref{data_example} we exemplify the application of our method to the analysis of time spent in 
intensive care units and weight gain. 
In Section \ref{discussion} we conclude with some final remarks. Technical proofs and derivations can be found in the Appendix.

\section{Compound expectation}\label{comp_exp}
  Let $Y$ be a random variable of interest with cumulative distribution function  
$F (y)$ and $\E(Y)<\infty$. The expectation of $Y$, $\mu$, can be expressed in terms 
of its quantile function $Q(p)$ as
\begin{equation*}
\mu=\E(Y)=\int_{-\infty}^{\infty} y \mathrm{d}F (y)=\int_0^1 Q(p) \mathrm{d}p.
\end{equation*} 
Given a set of specified proportions $\{\lambda_0,\ldots,\lambda_K\}$, 
where $\lambda_0=0$, $\lambda_K=1$ and $\lambda_{k-1}<\lambda_k$, for every $k=1,\ldots,K$, 
we can think of $\mu $ as a sum of $K$ components $\mu_1,\ldots,\mu_K$, that is, 
\begin{equation}\label{eq:component}
\mu=\sum_{k=1}^K \mu_k, \quad \mu_k=\int_{\lambda_{k-1}}^{\lambda_k}Q(p) \mathrm{d}p.
\end{equation} 

Each component $\mu_k$ amounts to a specific proportion of the total expectation, which
can be viewed as the contribution of the $k$th fraction of $Y$ to $\mu$. These contributions quantify the influence each 
fraction has on the overall mean, and allow to identify where efforts should be placed in order to, for example, 
increase it or decrease it. 
We define the $k$th contribution $c_k$ as follows
\begin{equation*}
c_k=\frac{\mu_k^++\mu_k^-}{\sum_{k=1}^K (\mu_k^++\mu_k^-)}, 
\end{equation*}
where $\mu_k^+=\int_{\lambda_{k-1}}^{\lambda_k}Q^+(p) \mathrm{d}p$ 
and $\mu_k^-=\int_{\lambda_{k-1}}^{\lambda_k}Q^-(p) \mathrm{d}p$, with $Q^+(p)=\max(Q(p),0)$ 
and $Q^-(p)=-\min (Q(p),0)$ the positive and negative parts of $Q(p)$ respectively. 
Components render a picture of the distribution of $Y$ in terms of its expected value, 
providing a measure of how much the specified fractions contribute to $\mu$. 

A measure of interest that can be easily derived from 
the decomposition presented in (\ref{eq:component}), is the expected value of $Y$ over each considered fraction. 
We define the compound expectation of $Y$ as the set 
$\{\overline{\mu}_k\}_{k=1}^K$, where for every $k=1,\ldots, K$,
\begin{equation*}
\overline{\mu}_k=\dfrac{\mu_k}{\lambda_{k}-\lambda_{k-1}}, \quad 
\mu=\sum_{k=1}^K (\lambda_{k}-\lambda_{k-1})\overline{\mu}_k,
\end{equation*}
and $\{\overline{\mu}_k\}_{k=1}^K$ measures the mean of $Y$ for every $k$th fraction.

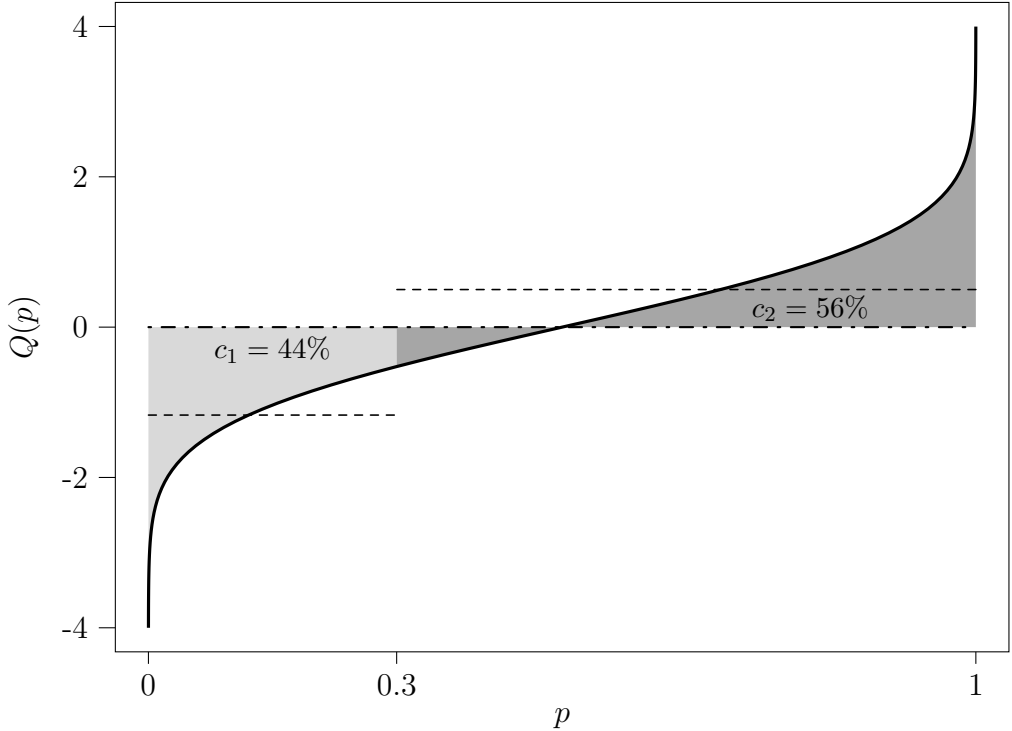
\begin{figure}[h!]
 \centering
 \begin{tikzpicture}[xscale=11,
      declare function={erf(\x)=%
      (1+(e^(-(\x*\x))*(-265.057+abs(\x)*(-135.065+abs(\x)%
      *(-59.646+(-6.84727-0.777889*abs(\x))*abs(\x)))))%
      /(3.05259+abs(\x))^5)*(\x>0?1:-1);}]
  \draw (-0.04,-4.32) -- (1.04, -4.32) -- (1.04, 4.32) -- (-0.04, 4.32) -- cycle;
  \foreach \y in {-4,-2,...,4}{
    \draw (-0.04, \y) -- (-0.06, \y) node[anchor=east]{\y};
    }
  \foreach \x/\z in {0/0,0.3/0.3,1/1}{
    \draw (\x, -4.32) -- (\x, -4.48) node[anchor=north]{\z};
    }
   \fill[gray!30] [domain=-4:-0.5244,variable=\y,samples=10] plot ({.5+.5*erf(\y/sqrt(2))},{\y}) -- 
		      (0.3,0) -- (0,0) -- cycle;
   \fill[gray!70] [domain=-0.5244:0,variable=\y,samples=10] plot ({.5+.5*erf(\y/sqrt(2))},{\y}) --
		      (0.3,0) -- cycle;
   \fill[gray!70] [domain=0:4,variable=\y,samples=40] plot ({.5+.5*erf(\y/sqrt(2))},{\y}) --
		      (1,0) -- cycle;
   \draw[very thick,domain=-4:4,variable=\y,samples=100] plot ({.5+.5*erf(\y/sqrt(2))},{\y});
   \draw[line cap=round, thick, dash pattern=on 1pt off 4pt on 5pt off 4pt] (0,0) -- (1,0);
   \draw[line cap=round, semithick,dashed] (0,-1.17) -- (0.3,-1.17);
   \draw[line cap=round, semithick,dashed] (1,0.5) -- (0.3,0.5);
   \node at (0.15, -0.3) {\small $c_1=44\%$};
   \node at (0.8, 0.25) {\small $c_2=56\%$};
   \node at (0.5, -5.2){$p$};
   \node [rotate=90] at (-0.15,0){$Q(p)$};
\end{tikzpicture}
  \caption{The light and dark gray shaded areas refer to lower and upper components, respectively,
  the dashed lines depict the compound expectation $\{\overline{\mu}_1,\overline{\mu}_2\}$, and the dot-dashed line depicts the total expectation. The contributions  
  are shown in text within each component's area.}
  \label{fig:examples}
\end{figure}

The components, contributions and compound expectation of an illustrative example are shown in Figure \ref{fig:examples}. 
We consider $Y$ that follows a standard normal distribution 
and we set $\{\lambda_0,\lambda_1,\lambda_2\}=\{0, 0.3, 1\}$. If we denote the standard normal 
quantile function as $\Phi^{-1}$, 
\begin{equation*}
\mu=\mu_1+\mu_2=\int_{0}^{0.3} \Phi^{-1}(p) \mathrm{d}p+\int_{0.3}^{1} 
\Phi^{-1}(p) \mathrm{d}p=-0.35+0.35=0
\end{equation*}
with $\{c_1,c_2 \}=\{0.44,0.56\}$ and $\{\overline{\mu}_1,\overline{\mu}_2\}=\{-1.17,0.5\}$. In this case, $30\%$ of $Y$ 
contributes $44\%$ to its total expectation $\mu$ and has a mean value of $-1.17$, while the remaining $70\%$ contributes only 
$56\%$ and has a mean value of $0.5$.

The number $K$ of components as well as the width between elements of the grid $\{\lambda_0,\ldots,\lambda_K\}$ 
determine the information one will retrieve from the underlying distibution. 
Note that the minimum grid $\{\lambda_0,\lambda_1\}$, i.e. $K=1$,  
corresponds to one component representing the total expectation $\mu$, 
whereas in the opposite end when $K\rightarrow \infty$, components become single quantiles.

In the following simple example, we highlight the advantages of the compound expectation over 
reporting single quantiles. 
Suppose we have two different groups of students, namely, Class A and Class B, that are graded according to a 
score scale that ranges from $0$ to $40$, and where an average of $20$ is considered the minimum satisfactory grade. 
While Class A has a total average score of $20.3$, 
Class B reaches only $19.6$, which is below the minimum desired. 

Aiming to improve the average score and performance in Class B, 
we take a further look at the scores' distributions in both groups. Table \ref{tab:example1} shows the compound 
expectation given a $K=10$ equally-sized components grid, 
$\{\lambda_0,\lambda_1,\ldots,\lambda_{10}\}=\{0,0.1,\ldots,1\}$, their corresponding 
contributions, and deciles of scores from Class A and Class B. While the deciles from both groups 
are extremely similar and do not show any apparent differences between the two distributions, 
the compound expectation clearly indicates that the worst $10\%$ 
of students in Class B have a lower average score than in Class A, identifying which students from Class B are 
more in need of an intervention (compared to Class A in this case). The bottom component summarizes all information contained in the lower tail of the 
distribution without the need of examining additional quantiles. This information is also conveyed by the 
contributions that are similar in both groups except for the bottom component, whose contribution in Class A is 
$2$ percentage points higher than in Class B. 

\begin{table}
  \caption{Compound expectation by deciles, overall means (last column), contributions (in $\%$) and deciles of final scores for Class A and Class B.}
\label{tab:example1}
  \begin{center}
    \begin{tabular}{ccccccccccc|c}
      Class A&k=1&k=2&k=3&k=4&k=5&k=6&k=7&k=8&k=9&k=10&$\mu$\tabularnewline
      \hline
      C. Expec.&10.9 &12.9 &15.3& 17.2& 18.9 &20.8 &22.5 &24.5 &27.4 &32.3&20.3\tabularnewline
      Contribut.&5.38&  6.36& 7.55&  8.48&  9.32& 10.3& 11.1&12.1 &13.5 &16.0&\tabularnewline
      Deciles&11.6& 13.7& 16.6 &18.0& 19.7& 21.7 &23.4 &25.8 &29.3&37.4&\tabularnewline
      \hline
      \tabularnewline
      Class B&k=1&k=2&k=3&k=4&k=5&k=6&k=7&k=8&k=9&k=10&$\mu$\tabularnewline
      \hline
      C. Expec.&6.10& 12.8 &15.1& 16.8 &18.7 &20.5 &22.5& 24.6 &27.3 &32.1&19.6\tabularnewline
      Contribut.& 3.10& 6.51 & 7.68  &8.55 & 9.52 &10.4& 11.5& 12.5 &13.9 &16.3&\tabularnewline
      Deciles&11.4 &13.6 &16.2&17.8 &19.5& 21.9& 23.3& 25.9 &29.5&38.5&\tabularnewline
      \hline
    \end{tabular}
  \end{center}
\end{table}
\section{Compound expectation in a regression framework}\label{reg_framework}
  Suppose that $x$ is an $m$-dimensional vector of covariates and we want to assess 
its effect on each component $\mu_k$. In this case, the components can be computed by means of 
the conditional quantile function, 
$Q(p|x)$, and for each value of $x$ we have $\mu(x)=\sum_{k=1}^K\mu_k(x)$, 
where $\mu(x)$ is the conditional expectation of $Y$ given $x$, 
and $\mu_k(x)=\int_{\lambda_{k-1}}^{\lambda_k} Q(p|x) \mathrm{d}p$ its $k$th conditional component. 
We focus hereafter on a 
specific class of all supposable $Q(p|x)$.

Let us assume that $Q(p|x)$ can be modeled in terms of 
$x$ and $p$ as follows
\begin{equation}\label{eq:factor}
Q(p|x)=t(x)^Tb(p),
\end{equation}
where $t:\mathbb{R}^m \rightarrow \mathbb{R}^q$ is a suitable transformation of $x$ and 
$b:(0,1) \rightarrow \mathbb{R}^q$ is coordinate-wise integrable, that is, its $i$th entry 
$b^i\in L^1(0,1)$ for $i=1,\ldots,q$. Models of this form have been extensively discussed
in the literature, e.g., in \cite{Koenker1978}, \cite{Efron1991}, \cite{Kim2007}, \cite{Cai2008}, 
\cite{Frumento2015} or \cite{Yuan2016}, among others.

Combining (\ref{eq:component}) and 
(\ref{eq:factor}) we have
\begin{align}\label{eq:regression_model}
\mu_k(x)&=\int_{\lambda_{k-1}}^{\lambda_k} t(x)^T b(p)\ \mathrm{d}p
=\int_{\lambda_{k-1}}^{\lambda_k} \left(\sum_{i=1}^q t^i(x)b^i(p)\right) \mathrm{d}p\nonumber\\
&= \sum_{i=1}^q \left(\int_{\lambda_{k-1}}^{\lambda_k} t^i(x)b^i(p)\ \mathrm{d}p \right)
= \sum_{i=1}^q \left(t^i(x)\int_{\lambda_{k-1}}^{\lambda_k} b^i(p)\ \mathrm{d}p\right)\nonumber\\
&=\sum_{i=1}^q t^i(x)B_{k}^i
=t(x)^T B_k,
\end{align}
where $B_k=(B_{k}^1,\ldots,B_{k}^q)^T$ with 
$B_{k}^i=\int_{\lambda_{k-1}}^{\lambda_k} b^i(p)\ \mathrm{d}p$, for $i=1,\ldots, q$.

Following the definition given in Section \ref{comp_exp}, we refer to 
$\{\overline{\mu}_k(x)\}_{k=1}^K$ as the conditional compound expectation (CCE) of $Y$, 
where 
\begin{equation*}
 \overline{\mu}_k(x)=t(x)^T \overline{B}_k, \quad
\overline{B}_k=\dfrac{B_k}{\lambda_{k}-\lambda_{k-1}},\text { for } k=1,\ldots, K.
\end{equation*}

Elements of $\{\overline{\mu}_k(x)\}_{k=1}^K$ are each presented as a regression model, where $\overline{B}_k$ measures average 
differences across values of $t(x)$ for a specific $k$th fraction of population,
making the class in (\ref{eq:factor}) of utmost convenience in our setting.

In accord with the remark in Section \ref{comp_exp}, the CCE generalizes two regression models of interest. For $K=1$, 
$$
\overline{\mu}(x)=\int_{0}^{1} t(x)^T b(p)\ \mathrm{d}p
= \sum_{i=1}^q \left(t^i(x) \int_{0}^{1} b^i(p)\ \mathrm{d}p \right)=\sum_{i=1}^q t^i(x)B^i,
$$
corresponds to linear regression, while 
$$
\lim_{\epsilon\to 0} \overline{\mu}^{(\epsilon)}_{k}(x)=\lim_{\epsilon\to 0}\frac{1}{\epsilon}\int_{\lambda_{k}-\epsilon}^{\lambda_{k}} t(x)^T b(p)\ \mathrm{d}p
=  \sum_{i=1}^q  \left(t^i(x) \lim_{\epsilon\to 0} \frac{1}{\epsilon} \int_{\lambda_{k}-\epsilon}^{\lambda_{k}} b^i(p)\ \mathrm{d}p \right)=\sum_{i=1}^q t^i(x)b^i(\lambda_{k}),
$$
corresponds to the $\lambda_k$-th regression quantile.

\section{Estimation and large sample properties}\label{asymp_properties}
  Let $Y_1,\ldots, Y_n$ be a set of independent replicates of a random variable $Y$ 
and $x_1,\ldots, x_n$ the corresponding $m$-dimensional vectors of fixed given covariates. We define 
$S=[L^1(0,1)]^q$, the $q$-dimensional $L^1(0,1)$ product space 
with norm $\lVert b \rVert_S = \sum_{i=1}^q \lVert b^i \rVert_{L^1}$, for all $b \in S$.

Based on the sample $\{(Y_j,x_j)\}_{j=1}^n$, we consider an estimator of the class of conditional quantile 
functions presented in (\ref{eq:factor}), $\widehat{Q}(p|x)=t(x)^T \hat{b}(p)$, with  
$\hat{b}\in S$. Our proposed estimator for 
$\mu_k(x)$ is then derived from (\ref{eq:regression_model}) and $\widehat{Q}(p|x)$, that is,
\begin{equation}\label{eq:cond_comp_estimator}
\hat{\mu}_k(x)=\int_{\lambda_{k-1}}^{\lambda_k}\widehat{Q}(p|x) \mathrm{d}p
=t(x)^T \hat{B}_k,
\end{equation}
where $\hat{B}_{k}=\int_{\lambda_{k-1}}^{\lambda_k}\hat{b}(p)\mathrm{d}p$ is an estimator of 
$B_{k}$, and all integrals are computed coordinate-wise. 

In this section we show that properties of $\hat{B}_k$ 
are naturally inherited from those of $\hat{b}$, and thus the choice of $\hat{b}$ 
will completely determine the features of our estimator. Proofs for Propositions~\ref{prop:unbiased}--\ref{prop:var_bound} 
can be found in Appendix A. 
Our first proposition characterizes the bias of $\hat{B}_k$ with respect to that of $\hat{b}$. 
\begin{proposition}\label{prop:unbiased}
If $\hat{b}$ is an unbiased estimator of $b$ almost everywhere, then
 $\hat{B}_k$ is an unbiased estimator 
 of $B_k$ for $k=1,\ldots, K$.
\end{proposition}
\begin{remark}
 For $\hat{B}_k$ to be unbiased, $\hat{b}$ needs not be unbiased almost everywhere 
but only the milder condition 
 $\int_{\lambda_{k-1}}^{\lambda_k} \left[ \E\{\hat{b}(p) -b(p) \} \right] \mathrm{d}p=0$ for  $k=1,\ldots, K$ 
 is necessary.
\end{remark}
Consistency and asymptotic distributional properties of $\hat{B}_k$ 
can also be derived from those of $\hat{b}$, as we present in the following two propositions.
\begin{proposition}\label{prop:consistency}
If $\hat{b}\in S$ is a consistent estimator of $b\in S$, then
 $\hat{B}_k$ is a consistent estimator of $B_k$ for $k=1,\ldots, K$.
\end{proposition}
\begin{proposition}\label{prop:asymptotic}
Let the sequence $(\tilde{b}_n)_{n\in\natus}\subset S$, where the $\sqrt{n}\tilde{b}_n=\sqrt{n}(\hat{b}_n-b)$ converge weakly to a 
$q$-dimensional Gaussian process $\{G(p)\}_{p\in(0,1)}$ with mean $m(p)=0$ for each $p\in (0,1)$, and 
matrix-valued covariance function 
$R(p,s)$, such that $G\in S$. If 
$\hat{B}_{k,n}=\int_{\lambda_{k-1}}^{\lambda_k}\hat{b}_n(p)\mathrm{d}p$, then $\sqrt{n}(\hat{B}_{k,n}-B_k)$
is asymptotically multivariate normal with mean 0 and covariance matrix 
$\Sigma_k=\int_{\lambda_{k-1}}^{\lambda_k}\int_{\lambda_{k-1}}^{\lambda_k} R(p,s) \mathrm{d}p\mathrm{d}s$, 
with $k=1,\ldots, K$.
\end{proposition}

Finally, the next proposition presents an upper bound for the variance of $\hat{B}^i_k$ 
in terms of that of $\hat{b}^i$.
\begin{proposition}\label{prop:var_bound}
Let $\hat{b}^i\in L^1(0,1)$ be a second order process, that is, $\E\{(\hat{b}^i)^2\}<\infty$ for each $p\in(0,1)$, then for $k=1,\ldots, K$
$$
\Var(\hat{B}^i_{k})
\leq(\lambda_k-\lambda_{k-1})\int_{\lambda_{k-1}}^{\lambda_k} \Var\{\hat{b}^i(p)\}\ \mathrm{d}p \text{ for } i=1,\ldots, q.
$$
\end{proposition}
\section{Grid of proportions}\label{grid}
  Both the width between elements of the grid of proportions $\{\lambda_0,\ldots,\lambda_K\}$ and its number $K$ of components 
define the information we aim to summarize from the underlying distribution of the data. For example, a grid with larger $K$, i.e., 
one with more components, 
may provide 
a more detailed summary, where more fractions are distinctly described. 
In practice, however, increasing the number of components might worsen estimation precision and could eventually 
provide less detailed (or uncertain) information. The reason behind this lies on the result presented in Proposition~\ref{prop:var_bound}, 
which upper bounds the variance of the CCE estimates $\hat{\mu}_k(x)$ by  
the variance of the underlying $\widehat{Q}(p|x)$. Segments on the variable's support where $\widehat{Q}(p|x)$ 
is more (less) precise will therefore yield better (worse) inferences to corresponding components.
Note that the bound in Proposition~\ref{prop:var_bound} is derived for each covariate separately (there referred to by the 
index $i$), and for every given grid, the estimates' precision might differ across covariates.

For illustrative purposes, we analyzed $300$ samples of a variable $Y$ generated from the heteroscedastic model 
$$
Y=2X+(3+Z)\varepsilon, 
$$
where $X\sim \chi^2(2)$, $Z\sim \text{Bernouilli} (0.5)$ and 
$\varepsilon\sim\normal (2,2)$.
CCEs were estimated given two different grids of proportions   
$$
   \Lambda^1=\{0, 0.1, 0.2, 0.3, 0.4, 0.5, 0.6, 0.7, 0.8, 0.9, 1\},\text{ and } \Lambda^2=\{0, 0.2, 0.4, 0.6, 0.8, 1\}
  $$
and using quantile regression (\cite{Koenker1978}) for the underlying conditional quantile function.
For every $k$th component $\mu_k(x)$, we considered a grid of order percentiles 
$\{p_0, p_1,\ldots, p_J\}$, where $p_j=\lambda_{k-1}+jh$ for $j=0,\ldots, J$, and $h=0.01$ the 
distance between consecutive elements of the grid. Note that 
$p_0=\lambda_{k-1}$ and $p_J=\lambda_{k-1}+Jh=\lambda_{k}$. 
A trapezoid rule approximation yielded
  $$
     \hat{B}^i_k=\int_{\lambda_{k-1}}^{\lambda_k}\hat{b}^i(p) \mathrm{d}p \approx 
     \sum_{j=0}^{J-1} \dfrac{\hat{b}^i(p_{j+1})-\hat{b}^i(p_{j})}{2}h,
  $$
and subsequently
\begin{equation*}
\hat{\mu}_k(x)=x^T\hat{B}_k= 
x^T \dfrac{h}{2}  \left(\hat{b}(p_{0})+2 \sum_{j=1}^{J-1}  \hat{b}(p_{j})+\hat{b}(p_{J})\right).
\end{equation*}

Confidence intervals for $\hat{\mu}_k(x)$ relied on the asymptotic normality results presented in \cite{Koenker1978} and 
on a discretized version of Proposition \ref{prop:asymptotic}. 
Exact expressions for $\Var(\hat{B}_k)$ and $\Var\{\hat{\mu}_k(x)\}$ for grid-based estimators can be found in Appendix~B.

\begin{figure}
\centering
  \begin{tikzpicture}[xscale=0.065, yscale=.15]
  \begin{scope}
    \DTLsetseparator{ }
    \DTLloaddb[noheader,keys={Component,Low,Mean,Up}]{X10}{./figures/x10_data.txt}
    
      \draw [->](-2,-7) -- (105,-7);
      \node at (50,-14) {Fractions of population (\%)};
      \foreach \x in {0,10,...,100}{
	\draw (\x, -7) -- (\x, -7.6) node[anchor=north]{\x};
      }
      
      \draw [->](0,-7) -- (0,26)  node [anchor=south] {(A)};
      \foreach \y in {-5,0,...,25}{
	\draw (0, \y) -- (-2, \y) node[anchor=east]{\y};
      }
    
      \def\interval{1}
      \def\sep{4}
      
    \foreach \comp/\shift in {0 - 0.1/5,0.1 - 0.2/15,0.2 - 0.3/25,0.3 - 0.4/35,0.4 - 0.5/45,0.5 - 0.6/55,0.6 - 0.7/65,0.7 - 0.8/75,0.8 - 0.9/85,0.9 - 1/95}{  
      \DTLforeach*[\DTLiseq{\c}{\comp}]{X10}{\c=Component,\l=Low,\m=Mean,\u=Up}{%
	\def\centralboxx{\shift}
	\draw[line cap=round,line join=round,thick] (-\interval+\centralboxx, \l) -- 
	(\interval+\centralboxx, \l) -- (\centralboxx, \l) 
	-- (\centralboxx, \u) -- (-\interval+\centralboxx, \u) -- (\interval+\centralboxx, \u); 
	\node at (\centralboxx, \m){\tiny$\bullet$};
      }
    }
  \end{scope}
  \begin{scope}[shift={(120,0)}]
    \DTLsetseparator{ }
    \DTLloaddb[noheader,keys={Component,Low,Mean,Up}]{Z10}{./figures/z10_data.txt}
    
      \draw [->](-2,-7) -- (105,-7);
      \node at (50,-14) {Fractions of population (\%)};
      \foreach \x in {0,10,...,100}{
	\draw (\x, -7) -- (\x, -7.6) node[anchor=north]{\x};
      }
      
      \draw [->](0,-7) -- (0,26) node [anchor=south] {(B)};
      \foreach \y in {-5,0,...,25}{
	\draw (0, \y) -- (-2, \y) node[anchor=east]{\y};
      }
    
      \def\interval{1}
      \def\sep{4}
      
    \foreach \comp/\shift in {0 - 0.1/5,0.1 - 0.2/15,0.2 - 0.3/25,0.3 - 0.4/35,0.4 - 0.5/45,0.5 - 0.6/55,0.6 - 0.7/65,0.7 - 0.8/75,0.8 - 0.9/85,0.9 - 1/95}{  
      \DTLforeach*[\DTLiseq{\c}{\comp}]{Z10}{\c=Component,\l=Low,\m=Mean,\u=Up}{%
	\def\centralboxx{\shift}
	\draw[line cap=round,line join=round,thick] (-\interval+\centralboxx, \l) -- 
	(\interval+\centralboxx, \l) -- (\centralboxx, \l) 
	-- (\centralboxx, \u) -- (-\interval+\centralboxx, \u) -- (\interval+\centralboxx, \u); 
	\node at (\centralboxx, \m){\tiny$\bullet$};
      }
    }
  \end{scope}
  \begin{scope}[shift={(0,-45)}]
    \DTLsetseparator{ }
    \DTLloaddb[noheader,keys={Component,Low,Mean,Up}]{X20}{./figures/x20_data.txt}
    
      \draw [->](-2,-7) -- (105,-7);
      \node at (50,-14) {Fractions of population (\%)};
      \foreach \x in {0,20,...,100}{
	\draw (\x, -7) -- (\x, -7.6) node[anchor=north]{\x};
      }
      
      \draw [->](0,-7) -- (0,26) node [anchor=south] {(C)};
      \foreach \y in {-5,0,...,25}{
	\draw (0, \y) -- (-2, \y) node[anchor=east]{\y};
      }
    
      \def\interval{1}
      \def\sep{4}
      
    \foreach \comp/\shift in {0 - 0.2/10,0.2 - 0.4/30,0.4 - 0.6/50,0.6 - 0.8/70,0.8 - 1/90}{  
      \DTLforeach*[\DTLiseq{\c}{\comp}]{X20}{\c=Component,\l=Low,\m=Mean,\u=Up}{%
	\def\centralboxx{\shift}
	\draw[line cap=round,line join=round,thick] (-\interval+\centralboxx, \l) -- 
	(\interval+\centralboxx, \l) -- (\centralboxx, \l) 
	-- (\centralboxx, \u) -- (-\interval+\centralboxx, \u) -- (\interval+\centralboxx, \u); 
	\node at (\centralboxx, \m){\tiny$\bullet$};
      }
    }
  \end{scope}
  \begin{scope}[shift={(120,-45)}]
    \DTLsetseparator{ }
    \DTLloaddb[noheader,keys={Component,Low,Mean,Up}]{Z20}{./figures/z20_data.txt}
    
      \draw [->](-2,-7) -- (105,-7);
      \node at (50,-14) {Fractions of population (\%)};
      \foreach \x in {0,20,...,100}{
	\draw (\x, -7) -- (\x, -7.6) node[anchor=north]{\x};
      }
      
      \draw [->](0,-7) -- (0,26) node [anchor=south] {(D)};
      \foreach \y in {-5,0,...,25}{
	\draw (0, \y) -- (-2, \y) node[anchor=east]{\y};
      }
    
      \def\interval{1}
      \def\sep{4}
      
    \foreach \comp/\shift in {0 - 0.2/10,0.2 - 0.4/30,0.4 - 0.6/50,0.6 - 0.8/70,0.8 - 1/90}{  
      \DTLforeach*[\DTLiseq{\c}{\comp}]{Z20}{\c=Component,\l=Low,\m=Mean,\u=Up}{%
	\def\centralboxx{\shift}
	\draw[line cap=round,line join=round,thick] (-\interval+\centralboxx, \l) -- 
	(\interval+\centralboxx, \l) -- (\centralboxx, \l) 
	-- (\centralboxx, \u) -- (-\interval+\centralboxx, \u) -- (\interval+\centralboxx, \u); 
	\node at (\centralboxx, \m){\tiny$\bullet$};
      }
    }
  \end{scope}
\end{tikzpicture}
  \caption{Estimated CCE (dots) with $95\%$ confidence intervals (line segments) 
  of $Y$ given $X=1$ (one unit increase) or $Z=1$ (binary variable category), given the grids $\Lambda_1=\{0, 0.1, 0.2, 0.3, 0.4, 0.5, 0.6, 0.7, 0.8, 0.9, 1\}$ 
  (panels A and B) and $\Lambda_2=\{0, 0.2, 0.4, 0.6, 0.8, 1\}$ (panels C and D).}
 \label{fig:gridssim}
\end{figure}
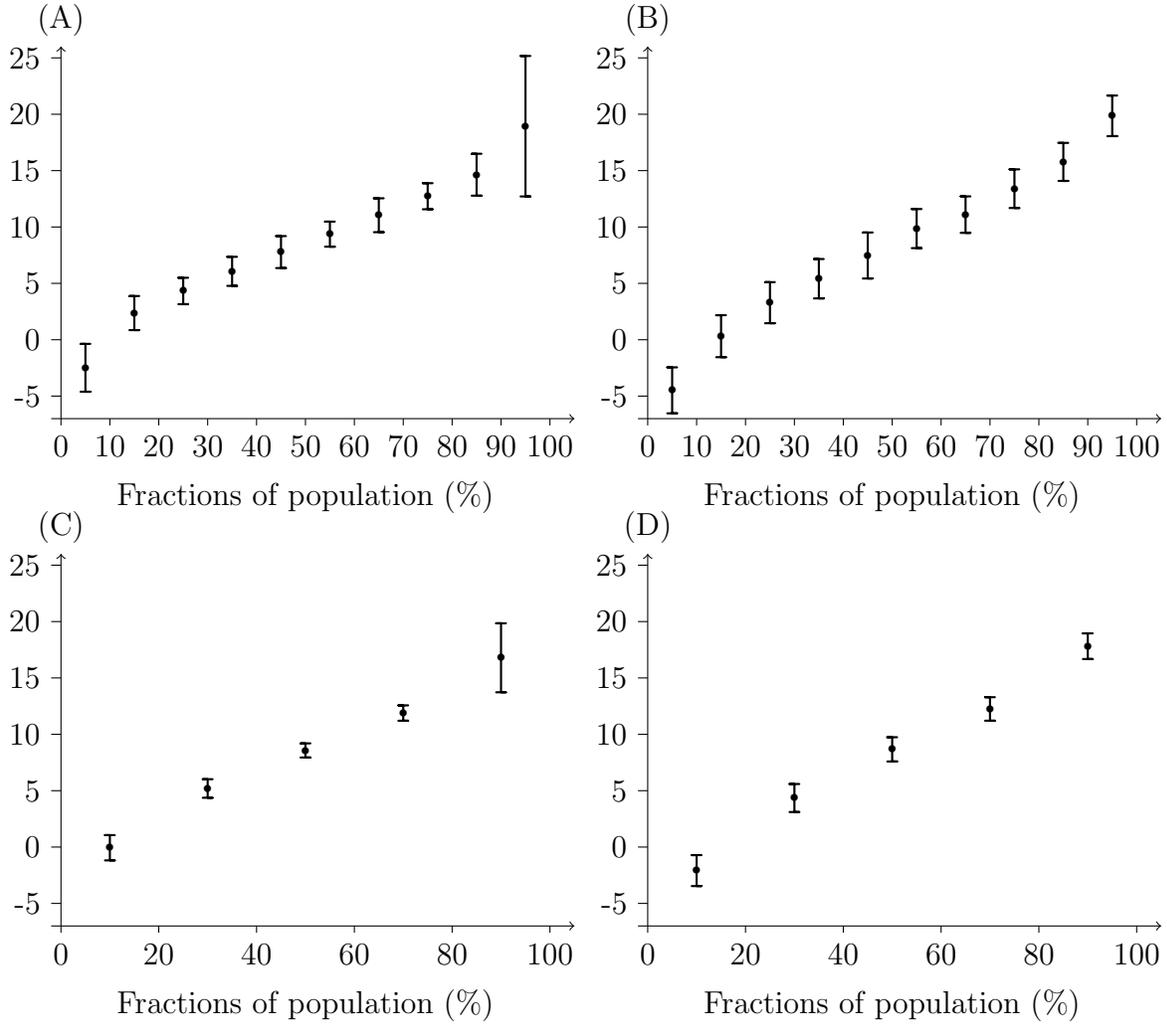

As we indicated before, $\Lambda^1$ might seem preferable when one aims to retrieve a more 
detailed summary of the data. However,  
the $95\%$ confidence interval for the top component of $Y$ given $X$ (with $Z$ set to $0$) is extremely large (Figure \ref{fig:gridssim} (A)), providing 
little information about the corresponding estimate. Intuitively, such imprecision relates to the skewness of the distribution of $X$, 
which translates to less precise estimates of higher order quantiles. 
If we consider results given $\Lambda^2$ instead (Figure \ref{fig:gridssim} (C)), estimation of the top component of $Y$ given $X$ 
is more precise, providing more detailed (although from a larger fraction) information.

We now observe that, given  $\Lambda^1$, estimation of the top component of $Y$ given $Z$ (with $X$ set to $0$) 
is considerably more precise (Figure \ref{fig:gridssim} (B)), specially when comparing it with other components' estimates. 
For the sake of completeness, we also show the compound expectation 
estimates of $Y$ given $Z$, given the coarser grid $\Lambda^2$ (Figure \ref{fig:gridssim} (D)).

We have summarized the effect that different grids might have on observed results, how it might vary 
across covariates, and underlined its close relation with both the chosen conditional quantile estimator and the underlying distribution of the data. 
Arguably, a finer grid of proportions ($\Lambda_1$) might not be optimal for obtaining a detailed summary of $Y$ given $X$ from our specific example, 
while it might be preferable when it comes to summarize $Y$ given $Z$. On the contrary, a coarser grid ($\Lambda_2$) would yield more
precise estimates but over 
larger fractions of population. Systematic approaches for obtaining optimal grids of proportions, in the sense of those that aim to find the best 
trade-off between precision and number of components, 
are highly dependable on the definition of how to exploit effectively the information contained in the data and are therefore 
not considered here. However, some comments on this point are included in our final remarks 
section (Section \ref{discussion}).  
\section{Real-data applications}\label{data_example}
  
\subsection{Average time spent in intensive care units}
The distribution of time spent in intensive care units (ICUs) is usually extremely skewed. A great majority 
of patients have relatively short stays whereas few present health complications and 
need to stay considerably longer. The latter are of interest because they experience worse medical outcomes and 
account for a major part of resources' allocation. 
In current studies, patients are often classified according to their length-of-stay and analyzed in separate groups,
(see for example, \cite{Stricker2003}, \cite{Trottier2007}, \cite{Kramer2010}), providing an ideal setting in which to apply our 
method. CCE estimates can thoroughly describe time spent in ICUs without the need of 
defining ad-hoc cut-points that hinder meaningful groups' comparisons, and will allow to identify which and 
how much patients influence average length-of-stay.

Given the grid $\Lambda=\{0, 0.1, 0.2, 0.3, 0.4, 0.5, 0.6, 0.7, 0.8, 0.9, 1\}$ and following the estimation strategy presented in Section \ref{grid}, 
we estimated the CCE of length-of-stay in 
ICUs of 31,828 middle-aged patients (45 to 65 years old) in Sweden, in relation with their gender, age and Charlson index 
(\cite{Charlson1987}). 
The Charlson index is an integer score scale that assesses severity of existing comorbidities
and is used as a prognostic indicator for disease outcome predictions. Higher Charlson scores 
indicate more severe comorbid conditions. Our study population was grouped in 2 
categories classified with respect to their Charlson scores, namely, those with Charlson index lower or equal to 1 (mild comorbidities) and those with Charlson index 
larger than 1 (severe comorbidities).
Further details on the data are available in \cite{rimes2015evolution}.

In Table \ref{tab:results} we show CCE estimates from where we could retrieve relevant 
information about average time spent in ICUs in relation with all considered factors.
For instance, we observed that among the first $60\%$ of the patients to leave the ICUs, 
average stay was of at most $0.85$ days for any age, gender and comorbid condition. 
However, the top component ranged between $6.76$ (for $45$ years old men with mild comorbidities) and $11.3$ 
(for $65$ years old women with severe comorbidities), yielding estimated total average stays ranging 
from $0.89$ to $1.5$ days. 
Interpreting these results in terms of contributions (shown in Table \ref{tab:contrib}), 
and for instance, in the reference group, 
we concluded that only 
$10\%$ of the patients (those who stayed the longest)  
accounted for $73.9\%$ of 
average time spent in ICUs. 

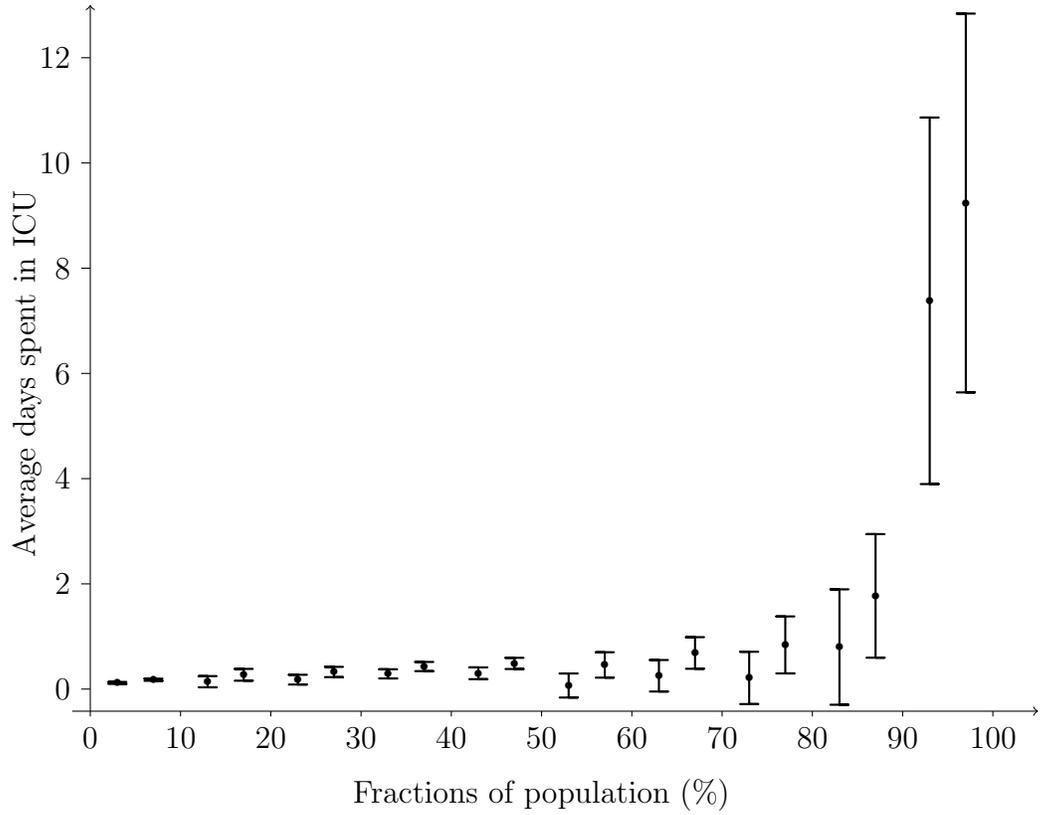
\begin{figure}
\centering
  \begin{tikzpicture}[xscale=0.12, yscale=0.7]
  \DTLsetseparator{ }
  \DTLloaddb[noheader,keys={Component,Group,Low,Mean,Up}]{CoMorb}{./figures/CoMorbplot_data.txt}
  
    \draw [->](-2,-0.42) -- (105,-0.42);
    \node at (50,-2) {Fractions of population (\%)};
    \foreach \x in {0,10,...,100}{
     \draw (\x, -0.42) -- (\x, -0.52) node[anchor=north]{\x};
    }
    
    \draw [->](0,-0.42) -- (0,13);
    \node at (-7, 6)[rotate=90]{Average days spent in ICU};
    \foreach \y in {0,2,...,12}{
      \draw (0, \y) -- (-1, \y) node[anchor=east]{\y};
    }
  
    \def\interval{1}
    \def\sep{4}
    
  \foreach \comp/\shift in {0.01 - 0.1/5,0.1 - 0.2/15,0.2 - 0.3/25,0.3 - 0.4/35,0.4 - 0.5/45,0.5 - 0.6/55,0.6 - 0.7/65,0.7 - 0.8/75,0.8 - 0.9/85,0.9 - 0.99/95}{  
    \DTLforeach*[\DTLiseq{\c}{\comp}]{CoMorb}{\c=Component,\g=Group,\l=Low,\m=Mean,\u=Up}{%
      \def\centralboxx{\shift+\g*\sep-1.5*\sep}
      \draw[line cap=round,line join=round,thick] (-\interval+\centralboxx, \l) -- 
      (\interval+\centralboxx, \l) -- (\centralboxx, \l) 
      -- (\centralboxx, \u) -- (-\interval+\centralboxx, \u) -- (\interval+\centralboxx, \u); 
      \node at (\centralboxx, \m){\tiny$\bullet$};
    }
  } 
  
\end{tikzpicture}
  \caption{Estimated CCE and $95\%$ confidence intervals of time 
spent in ICUs for 45 years old females across comorbid condition groups given the grid $\Lambda$. 
For each component we showed (from left 
to right) the mild and severe comorbidities groups.}
 \label{fig:grids}
\end{figure}

Average differences across comorbidity groups can also be retrieved and are summarized in Figure \ref{fig:grids}. 
We observed that on average, 
those with severe comorbidities experienced longer stays, a trend that increased consistently across components. 
The largest difference, estimated at $1.86$ days of average longer 
stay for those with severe comorbidities, was again observed amongst those on the top component, and represented a 
$35.2\%$ of the total average difference between the two groups, which was only $0.49$ days. 

Although estimates for the top component were not as precise as for other fractions of population, 
all results clearly indicate a high concentration amongst the $10\%$ with longest average stay in ICUs. Hence, the compound expectation 
and related quantities define a target group where to focus interventions that would yield a dramatic reduction on overall average 
stays and consequently average costs.

\begin{table}
  \caption{Estimated average and CCE given the grid $\Lambda$ of days spent in ICUs,
in relation to comorbid condition, age and gender. The reference group refers to 45 years old 
females with mild comorbidities. Standard errors are displayed in parenthesis.}
\label{tab:results}
  \begin{center}
    \begin{tabular}{c|cccc}
      \multicolumn{1}{c|}{\multirow{ 2}{*}{$\Lambda (\%)$}}&\multicolumn{1}{c}{\multirow{ 2}{*}{Reference}}&\multicolumn{1}{c}{Severe to}
      &\multicolumn{1}{c}{\multirow{ 2}{*}{Age}}&\multicolumn{1}{c}{Male to}\tabularnewline
      \multicolumn{1}{c|}{}&\multicolumn{1}{c}{}&\multicolumn{1}{c}{Mild comorb}
      &\multicolumn{1}{c}{}&\multicolumn{1}{c}{Female}\tabularnewline
      \hline
      $0-10$&$0.12(0.01)$&$0.06(0.01^{\star})$&$0.00(0.01^{\star})$&$0.00(0.01^{\star})$\tabularnewline
      $10-20$&$0.14(0.05)$&$0.13(0.01)$&$0.00(0.01^{\star})$&$-0.01(0.01)$\tabularnewline
      $20-30$&$0.18(0.05)$&$0.15(0.01)$&$0.00(0.01^{\star})$&$-0.01(0.01)$\tabularnewline
      $30-40$&$0.29(0.04)$&$0.14(0.01)$&$0.01(0.01^{\star})$&$0.00(0.01)$\tabularnewline
      $40-50$&$0.30(0.06)$&$0.19(0.01)$&$0.01(0.01^{\star})$&$0.00(0.01)$\tabularnewline
      $50-60$&$0.07(0.12)$&$0.39(0.03)$&$0.02(0.01^{\star})$&$0.02(0.02)$\tabularnewline
      $60-70$&$0.25(0.15)$&$0.43(0.04)$&$0.03(0.01^{\star})$&$0.02(0.03)$\tabularnewline
      $70-80$&$0.21(0.25)$&$0.63(0.07)$&$0.04(0.01^{\star})$&$-0.06(0.06)$\tabularnewline
      $80-90$&$0.80(0.56)$&$0.97(0.14)$&$0.07(0.01)$&$-0.20(0.13)$\tabularnewline
      $90-100$&$7.38(1.78)$&$1.86(0.41)$&$0.10(0.03)$&$-0.62(0.37)$\tabularnewline
      \hline
      Average&$0.97(0.17)$&$0.49(0.04)$& $0.03(0.01^{\star})$&$-0.08(0.04)$\tabularnewline
      \hline
      \multicolumn{5}{l}{\scriptsize{* denotes values that were smaller than the one displayed.}}
      \end{tabular}
      \end{center}
\end{table}

\begin{table}
  \caption{Contributions (in $\%$) of components of days spent in ICUs given the grid $\Lambda$,
in relation to comorbid condition, age and gender. The reference group refers to 45 years old 
females with mild comorbidities. Bootstrapped standard errors are displayed in parenthesis.}
\label{tab:contrib}
  \begin{center}
    \begin{tabular}{c|cccc}
      \multicolumn{1}{c|}{\multirow{ 2}{*}{$\Lambda (\%)$}}&\multicolumn{1}{c}{\multirow{ 2}{*}{Reference}}&\multicolumn{1}{c}{Severe to}
      &\multicolumn{1}{c}{\multirow{ 2}{*}{Age}}&\multicolumn{1}{c}{Male to}\tabularnewline
      \multicolumn{1}{c|}{}&\multicolumn{1}{c}{}&\multicolumn{1}{c}{Mild comorb}
      &\multicolumn{1}{c}{}&\multicolumn{1}{c}{Female}\tabularnewline
      \hline
      $0-10$&$1.17(0.29)$&$1.16(0.16)$&$0.07(0.07)$&$0.10(0.72)$\tabularnewline
      $10-20$&$1.56(0.55)$&$2.76(0.36)$&$1.20(0.34)$&$1.30(1.64)$\tabularnewline
      $20-30$&$2.00(0.54)$&$3.06(0.38)$&$2.28(0.40)$&$0.59(1.18)$\tabularnewline
      $30-40$&$3.23(0.61)$&$2.92(0.33)$&$2.74(0.41)$&$0.03(0.61)$\tabularnewline
      $40-50$&$3.34(0.66)$&$3.94(0.45)$&$3.76(0.61)$&$0.54(0.98)$\tabularnewline
      $50-60$&$0.76(0.64)$&$8.22(0.92)$&$6.99(0.99)$&$2.08(2.54)$\tabularnewline
      $60-70$&$2.82(1.31)$&$9.12(0.97)$&$9.23(1.06)$&$2.35(5.76)$\tabularnewline
      $70-80$&$2.36(1.71)$&$13.2(1.19)$&$15.2(1.39)$&$6.32(4.80)$\tabularnewline
      $80-90$&$8.89(4.49)$&$20.4(1.60)$&$25.7(1.87)$&$23.0(8.04)$\tabularnewline
      $90-100$&$73.9(5.64)$&$35.2(4.85)$&$32.8(4.94)$&$63.7(14.3)$\tabularnewline
      \hline
     \end{tabular}
  \end{center}
\end{table}

\subsection{Average weight gain amongst adolescents}
The obesity epidemic constitutes a major health burden in numerous countries, particularly in the United States. 
Health-care 
prevention strategies often focus on adolescents, so as to revert obesity at a young age and control the spread of the 
disease. In studies that aim to estimate the effect of a specific intervention, a commonly used measure is average 
weight gain over a certain period of time 
(see, for example, \cite{KornetvanderAa2017} for a review). 
Average weight gain estimates, however, are often misleading because they might easily take on negative values that will balance out with positive ones.
In these situations, the compound expectation proves a useful measure because it allows to distinguish between fractions of the population, clearly separating those 
who lost weight, or in other words, those with a negative weight gain, and quantifying their contribution to the total average. 

As an illustrative example, we estimated the CCE of weight gain amongst adolescents 
given the grid $\Lambda=\{0, 0.1, 0.2, 0.3, 0.4, 0.5, 0.6, 0.7, 0.8, 0.9, 1\}$, in relation with their age and obesity 
condition (obese versus non-obese). We analyzed a subset of $7,410$ adolescents from the National Longitudinal Survey of Youth 
(NLSY97) Cohort \citep{BureauofLaborStatistics2014} 
that followed young people in the United States (aged $13$ to $17$ on their first interview) 
through yearly or biennial interviews in which multiple aspects of their lives were reported. We considered their weight (in kilograms) at two different time points 
(as reported in 1997 and 2002) and computed their difference so that we could estimate average weight gain $5$ years after the 
first interview, and supposed they had been under a certain weight-control intervention that we aim to evaluate.

Table \ref{tab:results2} shows large disparities on estimated average weight gain across different fractions of 
population. For example in the reference group (non-obese $15$ years old), CCE estimates 
ranged from $-1.17$ (bottom component) to $33.3$ Kg (top component) with a total average of $13.3$ Kg of weight gain. 
Contributions of these components to the total average were $0.86\%$ and $24.6\%$, respectively (Table \ref{tab:contrib2}).

Regarding differences between the obese and non-obese groups we observe that, on average, the difference in 
weight gain was only of $1.27$ Kg, that is, hinting almost no effect of the intervention when comparing both groups. 
However, once divided into fractions of population, the bottom and top 
components represent a contribution of more than $20\%$ each to the total average difference, 
with estimates of $-6.08$ and $7.97$ Kg, respectively. In other words, while on average 
it might seem that the difference between the obese and non-obese groups is almost negligible, a closer inspection 
revealed  larger differences, and that those are highly concentrated around the $10\%$ of adolescents that lost and 
gained the most. In this case, the intervention would be producing almost opposite effects in the 
fractions corresponding to the bottom and top components. Those two groups, that represent the lower and upper tail of the 
distribution of weight gain difference, are comprehensively summarized by their respective components defining target groups 
in which to modify the intervention.

\begin{table}
  \caption{Estimated average and CCE of weight gain given the grid $\Lambda$,
in relation to obesity condition and age at baseline (1997). 
The reference group refers to 15 years old (mean age at baseline) non-obese adolescents. 
Standard errors are displayed in parenthesis.}
\label{tab:results2}
  \begin{center}
    \begin{tabular}{c|ccc}
      \multicolumn{1}{c|}{\multirow{ 2}{*}{$\Lambda (\%)$}}&\multicolumn{1}{c}{\multirow{ 2}{*}{Reference}}&
      \multicolumn{1}{c}{Obese to}
      &\multicolumn{1}{c}{\multirow{ 2}{*}{Age}}\tabularnewline
      \multicolumn{1}{c|}{}&\multicolumn{1}{c}{}&\multicolumn{1}{c}{Non-obese}
      &\multicolumn{1}{c}{}\tabularnewline
      \hline
      $0-10$&$-1.17(0.56)$&$-6.08(0.48)$&$-1.54(0.14)$\tabularnewline
      $10-20$&$3.96(0.33)$&$-1.6(0.29)$&$-1.78(0.11)$\tabularnewline
      $20-30$&$6.27(0.32)$&$-0.25(0.31)$&$-2.13(0.09)$\tabularnewline
      $30-40$&$8.49(0.40)$&$0.68(0.31)$&$-2.39(0.11)$\tabularnewline
      $40-50$&$10.7(0.42)$&$1.37(0.30)$&$-2.77(0.13)$\tabularnewline
      $50-60$&$13.1(0.44)$&$1.89(0.33)$&$-3.19(0.11)$\tabularnewline
      $60-70$&$15.9(0.44)$&$2.11(0.32)$&$-3.45(0.12)$\tabularnewline
      $70-80$&$19.2(0.48)$&$2.68(0.36)$&$-3.69(0.15)$\tabularnewline
      $80-90$&$23.3(0.69)$&$3.93(0.50)$&$-3.84(0.19)$\tabularnewline
      $90-100$&$33.3(0.91)$&$7.97(0.67)$&$-4.18(0.24)$\tabularnewline
      \hline
      Average&$13.3(0.17)$&$1.27(0.13)$& $-2.89(0.05)$\tabularnewline
      \hline
    \end{tabular}
  \end{center}
\end{table}

\begin{table}
  \caption{Contributions (in $\%$) of components of weight gain given the grid $\Lambda$,
in relation to to obesity condition and age at baseline (1997). The reference group 
refers to 15 years old (mean age at baseline) non-obese adolescents. Bootstrapped standard errors are displayed in parenthesis.}
\label{tab:contrib2}
  \begin{center}
    \begin{tabular}{c|ccc}
      \multicolumn{1}{c|}{\multirow{ 2}{*}{$\Lambda (\%)$}}&\multicolumn{1}{c}{\multirow{ 2}{*}{Reference}}&
      \multicolumn{1}{c}{Obese to}
      &\multicolumn{1}{c}{\multirow{ 2}{*}{Age}}\tabularnewline
      \multicolumn{1}{c|}{}&\multicolumn{1}{c}{}&\multicolumn{1}{c}{Non-obese}
      &\multicolumn{1}{c}{}\tabularnewline
      \hline
      $0-10$&$0.86(0.76)$&$21.3(2.52)$&$5.33(0.63)$\tabularnewline
      $10-20$&$2.93(0.22)$&$5.69(1.23)$&$6.14(0.33)$\tabularnewline
      $20-30$&$4.63(0.20)$&$0.86(0.78)$&$7.35(0.23)$\tabularnewline
      $30-40$&$6.27(0.27)$&$2.39(1.07)$&$8.24(0.26)$\tabularnewline
      $40-50$&$7.86(0.26)$&$4.78(0.97)$&$9.56(0.31)$\tabularnewline
      $50-60$&$9.70(0.30)$&$6.61(0.93)$&$11.0(0.25)$\tabularnewline
      $60-70$&$11.8(0.25)$&$7.38(0.78)$&$11.9(0.21)$\tabularnewline
      $70-80$&$14.2(0.28)$&$9.36(0.85)$&$12.8(0.28)$\tabularnewline
      $80-90$&$17.2(0.42)$&$13.8(1.33)$&$13.3(0.44)$\tabularnewline
      $90-100$&$24.6(0.64)$&$27.9(1.82)$&$14.4(0.62)$\tabularnewline
      \hline
     \end{tabular}
  \end{center}
\end{table}

\section{Final remarks}\label{discussion}
  A major strength of our approach lies on the possibility of choosing any quantile function estimator 
for the underlying distribution of the data, ensuring that all the compound expectation's estimation properties are to be inherited, 
and thus providing theoretically sound inference techniques. 
The large scope and variety of available methods, e.g., \cite{Portnoy2003} and \cite{Frumento2017} 
for censored data, or \cite{Koenker2004}, 
\cite{Geraci2007} and \cite{Wang2009a} for longitudinal data, among others,  
provide our method with flexibility to adapt to almost any scenario. 

Because quantiles are invariant to monotone transformations, 
models that estimate quantiles of these transformations from the original data might also be used. 
In such case, however, the estimated regression coefficients are often not directly interpretable and results are to 
be provided for each specific value of the covariates. 

The grid of proportions also provides flexibility to our measure, allowing to summarize the underlying distribution in as many 
fractions deemed appropriate. While in the present work we only considered a fixed given number and width of fractions of population in which to 
estimate the compound expectation, systematic approaches that 
choose optimally amongst all possible grids of proportions could also be explored. 
Systematic approaches could include, for example, finding the number and width of components that minimized 
the mean squared error between the underlying quantile function and the compound expectation estimates or that maintained 
the precision of the compound expectation estimates constant across components.

Finally, our measure has shown to be of great use in practical applications where average measures are of interest, and a 
more detailed look at the underlying distribution of the data reveals relevant information that might have been overlooked 
otherwise.
\section*{Appendix A}\label{app:proof}
\textbf{Proof of Proposition~\ref{prop:unbiased}.}
Because 
\begin{equation*}
\E(\hat{B}_k-B_k)=
\E\left\{\int_{\lambda_{k-1}}^{\lambda_k} \hat{b}(p) \mathrm{d}p-
\int_{\lambda_{k-1}}^{\lambda_k} b(p) \mathrm{d}p\right\}=
\int_{\lambda_{k-1}}^{\lambda_k} \left[ \E\{\hat{b}(p) \} -b(p) \right] \mathrm{d}p,
\end{equation*}
if $\E\{\hat{b}(p)-b(p)\}=0$ almost everywhere then $\E(\hat{B}_k-B_k)=0$, 
for $k=1,\ldots, K$.

\

\textbf{Proof of Proposition~\ref{prop:consistency}.}
Let $(\hat{b}_n)_{n\in\natus}\subset S$ be a sequence of consistent estimators of $b\in S$, that is, 
\begin{equation*}
 \lim_{n\to\infty} \Pr \left(\| \hat{b}_n-b \|_{S}>\varepsilon\right)= 0.
\end{equation*}

Consider the sequence $\hat{B}_{k,n}=\int_{\lambda_{k-1}}^{\lambda_k}\hat{b}_n(p)\mathrm{d}p$ of
estimators of $B_k$ for $k=1,\ldots, K$. 
Then, 
\begin{align*}
  \left\|\hat{B}_{k,n}-B_k\right\|_{\infty}&
  =\underset{i\in\{1,\ldots, q\}}{\max}\left|\int_{\lambda_{k-1}}^{\lambda_k}\left(\hat{b}_{n}^i(p) - b^i(p) \right)\mathrm{d}p\right|\\ 
  &\leq \sum_{i=1}^q \left|\int_{\lambda_{k-1}}^{\lambda_k} \left(\hat{b}_{n}^i(p) - b^i(p) \right)\mathrm{d}p\right|
  \leq \sum_{i=1}^q \int_{\lambda_{k-1}}^{\lambda_k} \left|\hat{b}_{n}^i(p) - b^i(p) \right|\mathrm{d}p\\
  &\leq \sum_{i=1}^q \int_{0}^{1} \left|\hat{b}_{n}^i(p) - b^i(p) \right|\mathrm{d}p
   = \left\| \hat{b}_n-b \right\|_{S},
\end{align*}
and because $\left\|\hat{B}_{k,n}-B_k\right\|_{\infty}\geq 0$, we have 
\begin{equation*}
 \lim_{n\to\infty} \Pr \left(\| \hat{B}_{k,n}-B_k \|_{\infty}>\varepsilon\right)= 0.
\end{equation*}
\

\textbf{Proof of Proposition~\ref{prop:asymptotic}.}
Consider the linear continuous operator 
  \begin{align*}
  T_k: S&\longrightarrow \reals^q\\
       b&\longmapsto \int_{\lambda_{k-1}}^{\lambda_k} b(p)\ \mathrm{d}p,
  \end{align*}
 where the integral is computed coordinate-wise. That $T_k$ is continuous follows easily from  
 \begin{align*}
  \lVert T_k(b)\rVert_{\infty}&=\left\lVert \int_{\lambda_{k-1}}^{\lambda_k} b(p)\ \mathrm{d}p\right\rVert_{\infty}
  =\underset{i\in\{1,\ldots, q\}}{\max} \left|\int_{\lambda_{k-1}}^{\lambda_k} b^i(p)\ \mathrm{d}p \right| 
  \leq \sum_{i=1}^q \lVert b^i \rVert_{L^1}= \lVert b\rVert_S.
 \end{align*}
  Because $\sqrt{n}\tilde{b}_n\rightarrow G$ in $S$, the continuous mapping theorem \cite[Th. 3.27, p. 54]{Kallenberg2006} yields 
 $$
 \sqrt{n}\int_{\lambda_{k-1}}^{\lambda_k} \tilde{b}_n(p)\ \mathrm{d}p\rightarrow
 \int_{\lambda_{k-1}}^{\lambda_k} G(p)\ \mathrm{d}p,\text{ that is, }
 \sqrt{n}(\hat{B}_{k,n}-B_k)\rightarrow\int_{\lambda_{k-1}}^{\lambda_k} G(p)\ \mathrm{d}p,
 $$
where $\rightarrow$ denotes weak convergence. 

Let us now compute the distribution of $\int_{\lambda_{k-1}}^{\lambda_k} G(p)\ \mathrm{d}p$. 
We have 
\begin{equation*}
  \int_{\lambda_{k-1}}^{\lambda_k} G(p)\ \mathrm{d}p=\lim_{L\to +\infty} \frac{\lambda_k-\lambda_{k-1}}{L}
\sum_{l=1}^LG_l,
\end{equation*} 
where $G_l=G(\lambda_{k-1}+lh)$ with $h=L^{-1}(\lambda_k-\lambda_{k-1})$.

Because $G=(G^1,\ldots, G^q)$ is a Gaussian process, we know that 
the vector 
$$
  \sum_{l=1}^LG_l=\left(\sum_{l=1}^LG^1_l,\ldots, \sum_{l=1}^LG^q_l\right)
$$
is multivariate normally distributed, with mean 0 and variance covariance matrix
\begin{equation*}
\Var\left(\sum_{l=1}^LG_l\right)
=\E\left(\sum_{l=1}^L\sum_{r=1}^LG_l^TG_r\right)
= \sum_{l=1}^L\sum_{r=1}^L\E (G_lG_r)=\sum_{l=1}^L\sum_{r=1}^L R(l,r).
\end{equation*}
The fact that $\sum_{l=1}^LG_l$ has mean $0$ follows immediately 
from $m(p)=0$ for all $p\in(0,1)$.

We consider now the normalised multivariate normal $Z_L\sim\mathcal{M}\normal(0,I_q)$ where $I_q$ is the 
$q$-dimensional identity matrix, 
$Z_L=\Sigma_{K,L}^{-1/2}L^{-1}(\lambda_k-\lambda_{k-1})\sum_{l=1}^LG_l$ 
and $\Sigma_{K,L}=L^{-2}(\lambda_k-\lambda_{k-1})^2\sum_{l=1}^L\sum_{r=1}^L R(l,r)$. It is clear that $Z_L\rightarrow Z$ with 
$Z\sim\mathcal{M}\normal(0,I_q)$. Applying now Slutsky's theorem we have
\begin{equation*}
 \frac{\lambda_k-\lambda_{k-1}}{L}\sum_{l=1}^LG_l\sim \mathcal{M}\normal(0,\Sigma_K), 
\end{equation*}
where 
$$
  \Sigma_k=\lim_{L\to +\infty}\left(\frac{\lambda_k-\lambda_{k-1}}{L}\right)^2\sum_{l=1}^L\sum_{r=1}^L R(l,r)=
  \int_{\lambda_{k-1}}^{\lambda_k}\int_{\lambda_{k-1}}^{\lambda_k} R(p,s) \mathrm{d}p\mathrm{d}s,
$$
concluding our proof.

\

\textbf{Proof of Proposition~\ref{prop:var_bound}.}
 \begin{align*}
  \Var(\hat{B}^i_{k})&=
  \E\left(\left|\hat{B}^i_{k}-\E(\hat{B}^i_{k})\right|^2\right)
 =\E\left(\left|\int_{\lambda_{k-1}}^{\lambda_k} \left[\hat{b}^i(p)-
 \E\{\hat{b}^i(p)\}\right] \mathrm{d}p\right|^2\right)\\
 & \leq \E\left(\left[\int_{\lambda_{k-1}}^{\lambda_k} \left|\hat{b}^i(p)-
 \E\{\hat{b}^i(p)\}\ \right|\mathrm{d}p\right]^2\right)\\
 &\stackrel{(1)}{\leq}\E\left[\int_{\lambda_{k-1}}^{\lambda_k} \left|\hat{b}^i(p)-
 \E\{\hat{b}^i(p)\}\right|^2\mathrm{d}p \int_{\lambda_{k-1}}^{\lambda_k} 1^2\mathrm{d}p\right]\\
 &=(\lambda_k-\lambda_{k-1})\int_{\lambda_{k-1}}^{\lambda_k} \E\left[\left|\hat{b}^i(p)-
 \E\{\hat{b}^i(p)\}\right|^2\right]\mathrm{d}p\\
 &=(\lambda_k-\lambda_{k-1})\int_{\lambda_{k-1}}^{\lambda_k}\Var\{\hat{b}^i(p)\}\mathrm{d}p,
 \end{align*}
where in (1) we use H\"{o}lders' inequality.
\section*{Appendix B}\label{app:inf}
  For any grid-based quantile function estimator like, for example, the one obtained from regression quantiles, 
the $q\times q$ variance-covariance matrices for the compound 
expectation regression coefficients $\hat{B}_{k}=(\hat{B}^1_{k},\ldots, \hat{B}^q_{k})$ for each $k=1,\ldots, K$ computed as in Section \ref{data_example} are given by
\begin{align*}
\Var\left( \hat{B}_{k}\right)&
=\Var\left\{\dfrac{h}{2}  \left(\hat{b}(p_{0})+2 \sum_{j=1}^{J-1}  
\hat{b}(p_{j})+\hat{b}(p_{J})\right)\right\}=
\dfrac{h^2}{4} \left\{\Var\left(\hat{b}(p_{0})\right)\right.\\
&+4\Var\left(\sum_{j=1}^{J-1}  \hat{b}(p_{j})\right)
+\Var\left(\hat{b}(p_{J})\right)
+2\Cov\left(\hat{b}(p_{0}),\sum_{j=1}^{J-1}  \hat{b}(p_{j})\right)\\
&+2\Cov^T\left(\hat{b}(p_{0}),\sum_{j=1}^{J-1}  \hat{b}(p_{j})\right)
+\Cov\left(\hat{b}(p_{0}),\hat{b}(p_{J})\right)
+\Cov^T\left(\hat{b}(p_{0}),\hat{b}(p_{J})\right)\\
&\left.+
2\Cov\left(\sum_{j=1}^{J-1}  \hat{b}(p_{j}),\hat{b}(p_{J})\right)+
2\Cov^T\left(\sum_{j=1}^{J-1}  \hat{b}(p_{j}),\hat{b}(p_{J})\right)\right\}\\
&=\dfrac{h^2}{4} \Bigg[ \Var\left(\hat{b}(p_{0})\right)+\Var\left(\hat{b}(p_{J})\right)+
\Cov\left(\hat{b}(p_{0}),\hat{b}(p_{J})\right)\\
&+\Cov^T\left(\hat{b}(p_{0}),\hat{b}(p_{J})\right)
+2\displaystyle\sum_{j=1}^{J-1} \left\{\Cov\left(\hat{b}(p_{0}),\hat{b}(p_{j})\right)+
\Cov^T\left(\hat{b}(p_{0}),\hat{b}(p_{j})\right)\right.\\
&+\left.\Cov\left(\hat{b}(p_{j}),\hat{b}(p_{J})\right)+
\Cov^T\left(\hat{b}(p_{j}),\hat{b}(p_{J})\right)\right\}
+4\displaystyle\sum_{j=1}^{J-1}\Var\left(\hat{b}(p_{j})\right)\\
&+4 \displaystyle\sum_{j=2}^{J-2} \sum_{l=j+1}^{J-1} \left\{ \Cov\left(\hat{b}(p_{l}),\hat{b}(p_{j})\right)+
\Cov^T\left(\hat{b}(p_{l}),\hat{b}(p_{j})\right) \right\} \Bigg],
\end{align*}
where $\Var\{\hat{b}(p_{j})\}$ denotes the variance-covariance matrix for the 
$p_j$th quantile regression coefficients 
and $\Cov\{\hat{b}(p_{j}),\hat{b}(p_{l})\}$ the 
covariance matrix between the $p_j$th and $p_l$th quantile regression coefficients.

The variance of each $k$th conditional component follows easily.
\begin{align*}
 \Var\{ \hat{\mu}_{k}(x)\}&=\Var\left(x^T\hat{B}_k\right)=
 \Var\left(\sum_{i=1}^q x^i\hat{B}_{k}^i\right)=\\
 &=\sum_{i=1}^q (x^i)^2\Var\left(\hat{B}_{k}^i\right)+
 2\sum_{i=1}^{q-1}\sum_{j=i+1}^q x^ix^j\Cov\left(\hat{B}_{k}^i,\hat{B}_{k}^j\right),
\end{align*}
where $\Var\left(\hat{B}_{k}^i\right)$ and $\Cov\left(\hat{B}_{k}^i,\hat{B}_{k}^j\right)$ are 
entries from 
the matrix $\Var\left( \hat{B}_{k}\right)$.

  \section*{Acknowledgements}\label{ack}
  The authors would like to thank Pol del Aguila Pla and Erin Gabriel for fruitful discussions that helped to improve the quality of this paper. 
  This project was partially funded by the KID doctoral grant from Karolinska Institutet.
\bibliography{CLE}
\end{document}